\RequirePackage{fix-cm}
\documentclass[twocolumn,envcountsect,envcountsame]{svjour3}          
\smartqed  
\usepackage{graphicx}
\usepackage[noadjust]{cite}

\usepackage{mathrsfs}
\usepackage{amsmath}
\usepackage{enumerate} 
\usepackage{amsxtra,latexsym, amscd} 
\usepackage{indentfirst}
\usepackage{amsfonts}
\usepackage{graphics}
\usepackage{multirow}
\usepackage{array}
\usepackage{verbatim}
\usepackage{mathabx}
\usepackage[usenames, dvipsnames]{xcolor}

\usepackage[hyperindex=true,
bookmarks=true,
pdftitle={}, pdfauthor={Abhishek Bhardwaj},
plainpages=false,
pdfpagelabels,
hyperfootnotes=false]{hyperref}
\hypersetup{
	colorlinks,
	linkcolor=blue,
	citecolor=OliveGreen,
	filecolor=magenta,      
	urlcolor=cyan,
}

\usepackage{nccmath}
		{\end{pmatrix}\end{medsize}}%
\newenvironment{mbmatrix}{\begin{medsize}\begin{bmatrix}}%
		{\end{bmatrix}\end{medsize}}%
\newcommand{\mc}{\mathcal}
\newcommand{\mf}{\mathfrak}
\newcommand{\mbb}{\mathbb}
\newcommand{\mbf}{\mathbf}

\newcommand{\norm}[1]{\left\lVert#1\right\rVert}

\newcommand{\set}[1]{\left\{#1\right\}}
\newcommand{\lr}[1]{\langle #1 \rangle}

\newcommand{\ul}[1]{\underline{#1}}

\newcommand{\LXR}{\lr{\ul{x}}}
\newcommand{\RX}{\mbb{R}\lr{\ul{x}}}
\newcommand{\SX}{\Sym\RX}
\newcommand{\mat}[1]{{#1}}
\newcommand{\dual}{\vee}
\newcommand{\invol}{\star}

\DeclareMathOperator{\Sym}{Sym}
\DeclareMathOperator{\Tr}{Tr}
\DeclareMathOperator{\tr}{tr}
\DeclareMathOperator{\CYC}{cyc}
\DeclareMathOperator{\s}{s}
\DeclareMathOperator{\rank}{rank}
\DeclareMathOperator{\I}{I}
\DeclareMathOperator{\II}{II}
\DeclareMathOperator{\st}{s.t.}

\def\M{\mathbf{M}}

\newcommand{\cyc}[2]{#1\overset{\CYC}{\sim} #2}

\usepackage{color}
\newcommand*{\mathcolor}{}
\def\mathcolor#1#{\mathcoloraux{#1}}
\newcommand*{\mathcoloraux}[3]{%
  \protect\leavevmode
  \begingroup
    \color#1{#2}#3%
  \endgroup
}
\newcommand{\sts}{\mathcolor{white}{\st}}

\begin{document}

\title{Noncommutative Polynomial Optimization
}


\author{Abhishek Bhardwaj \and Igor Klep \and Victor Magron }

\authorrunning{Bhardwaj, Klep \& Magron} 

\institute{A. Bhardwaj \at
				Australian National University\\
				\email{Abhishek.Bhardwaj@anu.edu.au}
           \and
           I. Klep \at
              University of Ljubljana \\
              \email{igor.klep@fmf.uni-lj.si}
           \and
           V. Magron \at
           		Laboratoire d'analyse et d'architectures des syst\`{e}mes \\
           		\email{Victor.Magron@laas.fr}           
}

\date{Received: date / Accepted: date}

\maketitle

\begin{abstract}
In this chapter we present the sums of Hermitian squares approach to noncommutative polynomial optimization problems. This is an extension of the sums of squares approach for polynomial optimization arising from real algebraic geometry. We provide a gentle introduction to the underlying theory of this methodology and highlight its importance.
\end{abstract}

\section{Introduction}
\label{intro}

Polynomial optimization problems (POP) are prevalent in many areas of modern science and engineering. The goal of POP is to minimize a given polynomial over a set defined by finitely many polynomial inequalities, a semialgebraic set. This problem is well known to be NP-hard \cite{laurent2009sums}, and has motivated research for more practical methods to obtain approximate solutions with high accuracy. Recent approaches to solving POPs use powerful representation results from real algebraic geometry, 
such as Putinar's Positivstellensatz \cite{putinar1993positive}. These methods aim to provide certificates of global positivity by reformulating a (positive) polynomial as a sum of other polynomials squared, acquiring a sum of squares (SOS) relaxation. Nowadays, Lasserre's Hierarchy \cite{lasserre2001global} is a standard method for handling SOS relaxations efficiently; it 
generates a sequence of bounds on the optimal values for POPs, with guaranteed convergence on a bounded constraint set. 

One can naturally extend the ideas of positivity and sums of squares to the noncommutative (nc) setting by replacing the commutative variables $z_{1}, \dotsc, z_{n}$ with noncommuting letters $x_{1}, \dotsc, x_{n}$. The extension to the noncommutative setting is an inevitable consequence of the many areas of science which regularly optimize functions with noncommuting variables, such as matrices or operators. For instance in control theory \cite{skelton1997unified, de2009engineering, camino2006solving}, matrix completion \cite{koltchinskii2011nuclear}, quantum information theory 
\cite{pal2009quantum, berta2016quantum, navascues2015characterizing, gribling2018bounds}, or quantum chemistry \cite{mazziotti2002variational, pironio2010convergent}.

The natural analogues of SOS representations in the noncommutative setting are {sums of Hermitian squares} (SOHS), which have similar representation results arising from real algebraic geometry \cite{helton2013free} and the tracial moment problems \cite{burgdorf2012truncated}. As in the commutative setting, these results have given rise to methods and tools (such as the software packages NCSOStools \cite{cafuta2011ncsostools}, and NCalgebra \cite{helton1996ncalgebra}) for the various noncommutative polynomial optimization problems (NCPOP), with many of the applications listed above readily benefiting from these tools. 

In this chapter we introduce the essentials of the sums of Hermitian squares approach to NCPOPs. We begin by introducing noncommutative analogues of the algebraic and geometric structures needed for optimization in \S \ref{S-Def}. We then formulate unconstrained NCPOPs in \S \ref{S-Unconstrained} for the eigenvalue minimization problem, and in \S \ref{S-trace} for the trace minimization problem. We also show how to represent and solve each of these problems as a semidefinite program, and \S \ref{S-Duality} considers the dual characterizations of these associated semidefinite programs. Including constraints in an NCPOP and the additional structures required to be able to obtain a solution are illustrated in \S \ref{S-Constraints}. In \S \ref{S-app} we present two important applications of NCPOPs and demonstrate how the SOHS methodology is applied. Finally in \S \ref{S-further-reading} we direct the reader to significant and recent research results and applications pertaining to NCPOPs.

\section{Definitions and preliminaries}\label{S-Def}

For a more detailed introduction to the concepts and ideas presented in this section, the reader is referred to \cite{burgdorf2016optimization}.

\subsection{Noncommutative polynomials}\label{SS-ncpoly}

Denote by $\ul{x}$ the noncommuting letters $x_{1}, \dotsc, x_{n}$. Let $\LXR = \lr{x_{1},\dotsc,x_{n}}$ be the \emph{free monoid} generated by $\ul{x}$, consisting of \emph{words} in $\ul{x}$, e.g., $x_{1}, x_{1}x_{2}, x_{2}x_{1}x_{1}x_{2} (=x_{2}x_{1}^{2}x_{2})$, etc. We will denote the empty word by $1$, and when we work with only two variables, we will use $x,y$ instead of $x_{1}, x_{2}$. Consider the free algebra $\RX$ of polynomials in $\ul{x}$ with coefficients in $\mbb{R}$. Its elements are called \emph{noncommutative \emph{(nc)} polynomials}, e.g., $t(x,y)=1+2 x+x^2+x y^2+2y^2+y^2 x+ yx^{2}y +y^4\in\mbb{R}\lr{x, y}$. Endow $\RX$ with the involution $f\mapsto f^{\invol}$ which fixes $\mbb{R}\cup\set{\ul{x}}$ pointwise, so that given a word $w=x_{i_1}\dots x_{i_r}$, $w^{\invol}$ is its reverse, i.e., $w^{\invol}=x_{i_r}\dots x_{i_1}$. This involution extends naturally to matrices of nc polynomials, i.e., given $\mat{V}$ with entries being nc polynomials, $(\mat{V}^{\invol})_{i,j} = (\mat{V}_{j,i})^{\invol}$.
The length of the longest word in a polynomial $f\in\RX$ is called the \emph{degree} of $f$ and is denoted $\deg(f)$, e.g., take $t$ as above, then $\deg(t)=4$. We write $\RX_{d}$ for all nc polynomials of degree at most $d$. The set of symmetric elements of $\RX$ is defined as $\SX = \set{f\in\RX : f^{\invol}=f}$. We employ the graded lexicographic ordering on all structures and objects we consider.

We write $\LXR_{d}$ for the set of all words in $\LXR$ of degree at most $d$, and we let $\mc{W}_{d}(\ul{x})=\mc{W}_{d}$  be the column vector of words in $\LXR_{d}$. The length of $\mc{W}_{d}$ is equal to  $\s(d,n) := \sum_{i=0}^{d} n^{i}$, which we write as $\s(d)$, when contextually appropriate. Given a polynomial $f\in\RX_{d}$, let $\mbf{f} = (f_{w})_{w\in \LXR_{d}}\in\mbb{R}^{\s(d)}$ be its vector of coefficients. It is clear that every polynomial $f\in\RX_{d}$ can be written as $f = \sum_{w\in \LXR_{d}}f_{w} w = \mbf{f}^{T}\mc{W}_{d} = \mc{W}_{d}^{\invol}\mbf{f}$. 

For $p\in\RX$, the nc polynomial $p^{\invol}p$ is called a \emph{Hermitian Square}. Given some $f\in\RX$, if there are other nc polynomials $g_{1}, \dotsc, g_{k}\in\RX$ such that we can write $f = g_{1}^{\invol}g_{1}+\dotsb+g_{k}^{\invol}g_{k}$, then $f$ is called a \emph{Sum of Hermitian Squares} (SOHS). The set of all such polynomials is 
\begin{equation}\label{eq-SOScone}
	\Sigma : = \set{ 
						\sum_{i=1}^{\ell} g_{i}^{\invol}g_{i} 
							\;\middle\vert\; 
						\ell\in\mbb{N}, g_{i}\in\RX 
					},
\end{equation}
and we write $\Sigma_{2d}$ for the elements of $\Sigma$ with degree less than or equal to $2d$, i.e., 
$\deg(g_{i})\leq d$ in the representation \eqref{eq-SOScone}.

\subsection{Positive Semidefinite Matrices}

Let $\mbb{S}_{r}$ denote the space of real symmetric matrices of size $r$; we will normally omit the subscript $r$ when we discuss matrices of arbitrary size, or if the size is clear from context. We will also denote by $\mat{I}_{r}$ the identity matrix of size $r$ (again omitting the subscript when appropriate).
Given $\mat{A}\in\mbb{S}$, $\mat{A}$ is positive semidefinite (psd) (resp. positive definite (pd)), if all eigenvalues of $\mat{A}$ are non-negative (resp. positive), and we write $\mat{A}\succeq0$ (resp. $\mat{A}\succ0$). Let $\mbb{S}^{\succeq}$ (resp. $\mbb{S}^{\succ}$) be the cone of psd (resp. pd) matrices. 
We denote by $\Tr(\mat{A})$ the trace $\sum_{i=1}^{r} A_{i,i}$ of the matrix $\mat{A} \in \mbb{S}_{r}$ and $\tr(\mat{A}) = \frac{1}{r} \Tr(\mat{A})$ is the
normalized trace. 
We write $\mat{\ul{A}} = (\mat{A}_{1}, \dotsc, \mat{A}_{n})\in\mbb{S}^{n}$, and given $q\in\mbb{R}\lr{\ul{x}}$, by $q(\mat{\ul{A}})$ we mean the evaluation of $q$ on $\mat{\ul{A}}$, i.e., replacement of the nc letters $x_{i}$ with the matrices $\mat{A}_{i}$ (note that considering $\ul{A}\in(\mbb{S}_{1})^{n}$ is the same as replacing the nc letters $x_{1}, \dotsc, x_{n}$ with commutative variables $z_{1}, \dotsc, z_{n}$). 

\subsection{Semidefinite Programming}\label{SS-SDP}

Semidefinite programs (SDP) are a well established tool for solving optimization problems, especially the ones arising from the ``sum of squares'' methodology. The theory of semidefinite programming is rich and substantial, and the interested reader is referred to \cite{wolkowicz2012handbook, anjos2011handbook}. We give the basic formulation of these programs below. 


Let $\lr{\cdot,\cdot}$ denote the scalar matrix product $\lr{\mat{A},\mat{B}} = \tr(\mat{B}^{T}\mat{A})$, then a semidefinite program has the standard primal form
\begin{equation}\label{eq-primalSDP}
    \begin{split}
        \xi =   & \sup_{X\in\mbb{S}_{r}} \lr{\mat{C},\mat{X}}, \\
                & \st \ \lr{\mat{A}_{j},\mat{X}} = b_{j},\ j=1,\dotsc,k,\\
                & \sts \ X \succeq 0,
    \end{split}
\end{equation}
with the matrix variable $\mat{X}$, and problem data $\mbf{b} = (b_{1},\dotsc,b_{k})^{T}\in\mbb{R}^{k}$, and $\mat{C}, \mat{A}_{1},\dotsc,\mat{A}_{k}\in\mbb{S}_{r}$. The dual problem to the above is the following \emph{dual SDP}
\begin{align*}
    \begin{split}
        \psi =  & \sup_{\mbf{y}\in\mbb{R}^{k}} \mbf{y}^{T}\mbf{b}, \\
		        & \st \ \mat{C}-\sum y_{i}\mat{A}_{i} \succeq 0.
    \end{split}
\end{align*}
When $\mat{C}$ is the zero matrix, (\ref{eq-primalSDP}) is called a feasibility program, and normally written as
\begin{align*}
    \begin{split}
		\sts \ & \mat{X}  \succeq 0, \\
		\st \  & \lr{\mat{A}_{j},\mat{X}} = b_{j},  \ j=1,\dotsc,k.
    \end{split}
\end{align*}

\section{Unconstrained Setting}\label{S-Unconstrained}

One of the standard applications of NCPOP is \emph{eigenvalue optimization}, 
\begin{equation}\label{eq-ev}
	\lambda_{\min}(f) := \inf \set{ 
									    \mbf{v}^{T}f(\mat{\ul{A}})\mbf{v} 
										    \;\middle\vert\; 
									    \mat{\ul{A}}\in\mbb{S}^{n}, \text{ and } \norm{\mbf{v}}=1 
									},
\end{equation}
which aims to find the minimum possible eigenvalue of an nc polynomial $f$ over $\mbb{S}^{n}$.
Reformulating, \eqref{eq-ev} is equivalent to 
\begin{equation}\label{eq-ev.sup}
	\lambda_{\min}(f) = \sup\set{ 
									\lambda 
										\;\middle\vert\; 
									\forall \mat{\ul{A}}\in\mbb{S}^{n}, f(\mat{\ul{A}}) -\lambda \mat{I}\succeq0  
								},
\end{equation}
and as the following theorem shows, for nc polynomials, matrix positivity (as required in \eqref{eq-ev.sup}) is equivalent to membership in $\Sigma$.

\begin{theorem}[\protect{\cite{mccullough2001factorization, helton2002positive, mccullough2005noncommutative}}]\label{T-matrixpositiveSOHS}
Given $f\in\RX$, $f$ is \emph{matrix positive}, i.e., $f(\mat{\ul{A}})\succeq0$ for all $\mat{\ul{A}}\in\mbb{S}^{n}$, if and only if $f$ is a sum of Hermitian squares, i.e., $f\in\Sigma$.
\end{theorem}
Using this theorem, we can frame eigenvalue minimization as the following SOHS program
\begin{equation}\label{eq-ev.sohs}
	\lambda_{\min}(f) = \sup\set{ 
									\lambda 
										\;\middle\vert\; 
									f - \lambda\in\Sigma
								},
\end{equation}
which can be solved through the use of semidefinite programming.

\subsection{SOHS decomposition via SDP}\label{SS-PrimalSDP}

Let us now consider how to test membership into $\Sigma$. We note that we can write every $f\in\SX_{2d}$ as $f = \mc{W}_{d}^{\invol}\mat{G}_{f}\mc{W}_{d}$, for some appropriate symmetric matrix $\mat{G}_{f}$, which is known as a \emph{Gram matrix} of $f$. As the example below shows, the Gram matrix is not unique.

\begin{example}
	Consider $t(x,y)=1+2 x+x^2+x y^2+2y^2+y^2 x+yx^2 y+y^4$. We can represent $t$ as
	$$
	t(x,y)=
	\mc{W}_{2}^{\invol}
	\begin{mbmatrix}
		1&	1    &	0    &	a &	0 &	0 &	b  \\
		1&	1-2a &	0    &	0 &	0 &	0 &	1-c\\
		0&	0    &	2-2b &	0 &	0 &	c &	0  \\
		a&	0    &	0    &	0 &	0 &	0 &	0  \\
		0&	0    &	0    &	0 &	1 &	0 &	0  \\
		0&	0    &	c    &	0 &	0 &	0 &	0  \\
		b&	1-c  &	0    &	0 &	0 &	0 &	1
	\end{mbmatrix}
	\mc{W}_{2}
	$$
	for any $a,b,c\in\mbb{R}$.
\end{example}

\begin{lemma}[\protect{\cite[Proposition 1.16]{burgdorf2016optimization}}]
\label{L-SOHSdecomp}
	Let us consider a polynomial $f\in\SX_{2d}$. Then $f\in\Sigma$ if and only if $f$ has a psd Gram matrix, i.e., 
	there exists a psd matrix $\mat{G}_{f}$ such that
	$$
	f = \mc{W}^{\invol}_{d}\mat{G}_{f}\mc{W}_{d}.
	$$
\end{lemma}

\newpage

\begin{example}
	Let $t(x,y)=1+2 x+x^2+x y^2+2y^2+y^2 x+ yx^2 y+ y^4$, with its Gram matrix 
	$$
	\mat{G}(a,b,c) = 
	\begin{mbmatrix}
		1&	1    &	0    &	a &	0 &	0 &	b  \\
		1&	1-2a &	0    &	0 &	0 &	0 &	1-c\\
		0&	0    &	2-2b &	0 &	0 &	c &	0  \\
		a&	0    &	0    &	0 &	0 &	0 &	0  \\
		0&	0    &	0    &	0 &	1 &	0 &	0  \\
		0&	0    &	c    &	0 &	0 &	0 &	0  \\
		b&	1-c  &	0    &	0 &	0 &	0 &	1
	\end{mbmatrix}
	$$ 
	It is readily seen that $G(a,b,c)\succeq0$ if and only if $a=b=c=0$. So,
	we know that $t$ is an SOHS, and indeed we can write
	$$
	t(x,y) = (1+x+y^2)^{\invol}(1+x+y^2) + (xy)^{\invol}(xy). 
	$$
\end{example}

Computing an SOHS decomposition for any $p\in\SX_{2d}$ amounts to comparing ``symmetrized coefficients'' (i.e. $p_{w}+p_{w^{\invol}}$) with coefficients of $\mc{W}_{d}^{\invol}\mat{G}\mc{W}_{d}$, for some psd $\mat{G}$. Such a psd $\mat{G}$ can be obtained as the solution of the following semidefinite feasibility program,
\begin{align*}
    \begin{split}
		\sts \ & G \succeq0,\\
		\st \    & \langle \mat{A}_{w}, G \rangle = p_{w} + p_{w}^{\invol}, \ w\in \mc{W}_{2d}    
    \end{split}
\end{align*}
where 
$$
(\mat{A}_{w})_{u, v} = 
\begin{cases}
	2, \text{ if } u^{\invol}v=w, \text{ and } w=w^{\invol}, \\
	1, \text{ if } u^{\invol}v\in\set{w, w^{\invol}} \text{ and } w\neq w^{\invol},\\
	0, \text{ otherwise}.
\end{cases}
$$
Combining Theorem \ref{T-matrixpositiveSOHS} with Lemma \ref{L-SOHSdecomp} gives the SDP form of the SOHS program \eqref{eq-ev.sohs},
\begin{align}\label{eq-SOHShierarchy}
	\begin{split}
		\lambda_{\min}(f) = & \sup \ \lambda,\\
		        & \st \ f-\lambda  = \mc{W}_{d}^{\invol}G\mc{W}_{d}, \\
		        & \sts \  G\succeq0.
	\end{split}
\end{align}


\section{Trace Optimization}\label{S-trace}

Trace optimization is an important branch of noncommutative optimization, with many different applications. Of particular importance is the quantum mechanical setting, where 
one investigates \emph{statistical averages} of the system, which can be naturally represented using the trace.

\subsection{Cyclic Equivalence}\label{SS-Cyclic equivalences}

Given polynomials $p,q\in\RX$, the \emph{commutator} of $p$ and $q$ is defined as $[p,q]:=pq-qp.$ It is well known that a trace zero matrix is a (sum of) commutator(s). So when optimizing over the trace, we remove commutators through the equivalence relation defined below.

\begin{definition}
	Two polynomials $f,g\in\RX$ are called \emph{cyclically equivalent} ($\cyc{f}{g}$) if $f-g$ is a sum of commutators, i.e.,
	$$
	f - g =\sum_{i=1}^{\ell} p_{i}q_{i}-q_{i}p_{i} \text{ for some } \ell\in\mbb{N} \text{ and } p_{i}, q_{i}\in\RX.
	$$
\end{definition}

Two words $u,v\in\LXR$ are cyclically equivalent if and only if there exists words $w_{1}, w_{2}\in\LXR$ such that $u=w_{1}w_{2}$ and $v=w_{2}w_{1}$. Using this with linearity of the commutator, we can check cyclic equivalence of nc polynomials, by comparing (the sum of) coefficients of cyclically equivalent words. This is formalized in the next proposition.

\begin{proposition}[\protect{\cite[Proposition 1.51]{burgdorf2016optimization}}]\label{P-cyclictolinear}
	Two nc polynomials $f = \sum_{w\in\LXR} f_{w} w$ and $g = \sum_{w\in\LXR} g_{w} w$ ($f_{w}, g_{w}\in\mbb{R}$) are cyclically equivalent if and only if for all $v\in\LXR$ one has
	$$ 
		\sum_{ \substack{ w\in\LXR \\ \cyc{w}{v} } } f_{w} = \sum_{ \substack{ w\in\LXR \\ \cyc{w}{v} }} g_{w} .
	$$ 
\end{proposition}

\subsection{Trace Minimization}
With cyclic equivalence, we now have the following equality
\begin{align*}
    \begin{split}
        &\inf\set{ 
	    		    \tr(f(\mat{\ul{A}})) 
		    	        \;\middle\vert\;
			        \mat{\ul{A}}\in\mbb{S}^{n}, f\in\RX 
		        }= \\
	    &\inf\set{
		    		\tr(g(\mat{\ul{A}})) 
			    	    \;\middle\vert\; 
				    \mat{\ul{A}}\in\mbb{S}^{n}, g\in\RX, \text{ and } \cyc{g}{f} 
			    }.
	\end{split}
\end{align*}
Hence, if we consider the trace minimization problem
\begin{align*}
    \begin{split}
    	\tr_{\min}(f) &:= \inf \set{ 
	    										\tr(f(\mat{\ul{A}})) 
		    									\;\middle\vert\; 
			    								\mat{\ul{A}}\in\mbb{S}^{n}
			    							} \\
	    &= \sup \set{ 
		    				\tau 
			    			\;\middle\vert\; 
				    		\tr(f(\mat{\ul{A}})) - \tau \geq 0, \forall \mat{\ul{A}}\in\mbb{S}^{n}
					}
	\end{split}
\end{align*}
to generate a relaxation as in \eqref{eq-ev.sohs}, we have to consider all nc polynomials, which are cyclically equivalent to something in $\Sigma$. So we define
\begin{equation*}
	\Theta := \set{ 
						f\in\RX 
							\;\middle\vert\; 
						\exists g\in\Sigma \text{ with } \cyc{f}{g} 
				    }
\end{equation*}
and obtain the relaxation
\begin{equation}\label{eq-trSOHS}
	\tr_{\Theta}(f) = \sup\set{ 
									\tau 
										\;\middle\vert\; 
									f -\tau\in\Theta  
								}.
\end{equation}

In contrast to eigenvalue optimization, we know only that $\tr_{\Theta}(f)\leq \tr_{\min}(f)$, with the inequality being strict in general.

\begin{example}
	Consider the nc Motzkin polynomial, 
	$$
	M_{\text{nc}}(x,y) = xy^{4}x + yx^{4}y - 3xy^{2}x +1\in \Sym\mbb{R}\lr{x,y}
	$$
	It is known that $M_{\text{nc}}$ is trace positive, i.e., for every $A,B\in\mbb{S}$, $\tr(M_{\text{nc}}(A,B))\geq 0$ (see \cite[Example 4.4]{klep2008connes}), from which we can see that $\tr_{\min}(M_{\text{nc}})=0$. However, $\tr_{\Theta}(M_{\text{nc}})=-\infty$, for if there exists a $\tau\in\mbb{R}$ with $M_{\text{nc}}-\tau\in\Theta$, then (considering $A,B\in\mbb{S}_{1}$) the (commutative) Motzkin polynomial would be a sum of squares, and this is well known to be false \cite[Proposition 1.2.2]{zbMATH05255035}.
\end{example}

\subsection{SDP formulation}

For an nc polynomial $f\in\RX_{2d}$, we can search for a cyclic SOHS decomposition as before with the following SDP
\begin{align*}
    \begin{split}
		\sts \ & G \succeq 0,\\
		\st \    & \cyc{f}{\mc{W}_{d}^{\invol}G\mc{W}_{d}},
    \end{split}
\end{align*}
where we point out that checking cyclic equivalence is done through a set of linear constraints, as detailed in Proposition \ref{P-cyclictolinear}.

With this we obtain the SDP form of \eqref{eq-trSOHS}, to approximate $\tr_{\min}(f)$,
\begin{align}\label{eq-trhierarchy}
	\begin{split}
		\tr_{\Theta}(f) = &  \sup \ \tau, \\
                          & \st \ \cyc{f-\tau}{\mc{W}_{d}^{\invol}G\mc{W}_{d}}, \\
                          & \sts \ G\succeq 0.
	\end{split}
\end{align}


Despite the existence of examples where $\tr_{\Theta}(f) < \tr_{\min}(f)$, for most problems $\tr_{\Theta}$ is more computationally viable and offers a good approximation to $\tr_{\min}$.

\begin{example}
	Consider $t(x,y)=1+2 x+x^2+x y^2+2y^2+y^2 x+yx^2 y+y^4$
	, which we know is an SOHS, and so necessarily, $\tr_{\min}(t) \geq 0$. Using \texttt{NCSOStools} we compute $\tr_{\Theta}^{(2)}(t)$ with the following commands
	\begin{verbatim}
		>> NCvars x y 
		>> t = 1+2*x+x^2+x*y^2+2*y^2+y^2*x+y*x^2*y+y^4;
		>> d = 2;
		>> opt_d = NCtraceOpt(t, {}, 2*d);
	\end{verbatim}
	and obtain $\tr_{\Theta}^{(2)}(t) \approx 6.967\times 10^{-10}$. Thus we conclude $\tr_{\min}(t) = 0$.
\end{example}

\section{Duality}\label{S-Duality}

As is known from the duality theory of SDP \cite{wolkowicz2012handbook}, the dual optimum provides a lower bound on the primal optimum (in the notation of \S \ref{SS-SDP}; $\psi\leq\xi$, see \cite[\S 3]{vandenberghe1996semidefinite} for an example). With this in mind, we now describe the dual formulations of \eqref{eq-ev.sohs} and \eqref{eq-trSOHS}.

\subsection{Eigenvalue Dual}

Let $\Sigma^{\dual}_{2d}$ be the dual cone of $\Sigma_{2d}$. This dual cone is characterized as 
\begin{equation*}
	\Sigma^{\dual}_{2d} : = \set{ 
									L:\RX_{2d}\rightarrow\mbb{R} 
										\;\middle\vert\; 
									\begin{aligned}		    		   L \text{ is linear}, \\ 
								    L(f) = L(f^{\invol}), \\
									\forall f\in\Sigma_{2d}, L(f)\geq0 
									\end{aligned}
								}
\end{equation*}
Following a standard Lagrangian duality argument (cf. \cite[Chapter 4]{wolkowicz2012handbook}) for eigenvalue optimization, we obtain the dual program to \eqref{eq-SOHShierarchy},
\begin{align}\label{eq-evdual}
	\begin{split}
		\Lambda_{\min}(f) = & \inf \ L(f), \\
		                    & \st \ L \in \Sigma^{\dual}_{2d}, \\
		                    & \sts \ L(1) = 1.
	\end{split}
\end{align}

Linear functionals $L$ on $\RX_{2d}$ which are \emph{symmetric} ($L(f)=L(f^{\invol})$), can be represented as \emph{nc Hankel matrices}.

\begin{definition}
	A matrix $\M \in\mbb{S}_{\s(d)}$, which is indexed by words $u,v\in\LXR_{d}$, is an \emph{nc Hankel matrix} if 
	$$
	\M_{u,v} = \M_{r,s}, \quad \text{whenever} \quad u^{\invol}v = r^{\invol}s.
	$$
\end{definition}

\begin{definition}
	Let $L:\RX_{2d}\rightarrow\mbb{R}$ be a symmetric linear functional. The \emph{nc Hankel matrix} $\M_d(L)$, associated to $L$, is defined as
	\begin{equation*}
		\M_d(L)_{u,v} = L(u^{\invol}v).
	\end{equation*}
\end{definition}
Let $p=\sum p_{w} w\in\RX_{d}$ and $q = \sum q_{w} w\in\RX_{d}$ be two nc polynomials. Then it is easy to see that
\begin{equation}\label{eq-hankelpositive}
L(p^{\invol}q) = \mbf{p^{T}} \M_d(L) \mbf{q},
\end{equation}
and in fact $\M_d(L)$ is the unique matrix with this property. Moreover, from \eqref{eq-hankelpositive}, we see that for all $p\in\RX$, $L(p^{\invol}p)\geq0$ if and only if $\M_d(L)\succeq0$. This gives us an alternative representation of the dual cone $\Sigma^{\dual}_{2d}$, and lets us rewrite \eqref{eq-evdual} as 
\begin{equation}\label{eq-evSDPdual}
	\begin{split}
		\Lambda_{\min}(f) = & \inf \ \lr{\M_d(L), \mat{G}_{f}}, \\
		& \st \ \M_{d}(L)\succeq0, \\
		& \sts \ \M_d(L)_{1,1} = 1, \\
		& \sts \ \M_d(L)_{u,v} = \M_d(L)_{r,s} \text{ for all } u^{\invol}v = r^{\invol}s,
	\end{split}
    \hspace{10000pt minus 1fil}
\end{equation}
where $\mat{G}_{f}$ is a Gram matrix of $f$.

\subsection{Trace Dual}

For trace optimization, we can similarly obtain the dual program to \eqref{eq-trhierarchy},
\begin{align}\label{eq-trdual}
	\begin{split}
		L_{\Theta}(f) = & \inf \ L(f), \\
		& \st \ L\in \Theta^{\dual}_{2d}, \\
		& \sts \ L(1)=1.
	\end{split}
\end{align}
The dual cone $\Theta^{\dual}_{2d}$ can be characterized similarly to $\Sigma^{\dual}_{2d}$. In fact, since $\Sigma_{2d}\subseteq\Theta_{2d}$, we know $\Theta^{\dual}_{2d}\subseteq\Sigma^{\dual}_{2d}$. So the dual cone $\Theta^{\dual}_{2d}$ will certainly consist of linear functionals $L:\RX_{2d}\rightarrow\mbb{R}$ which are symmetric, and $L(\Sigma_{2d})\subseteq[0,\infty)$. Moreover, for a functional to be well defined under cyclic equivalences, $\Theta^{\dual}_{2d}$ will consist of the \emph{tracial} linear functionals, which satisfy $L([p,q])=0$, for $p,q\in\RX_{d}$. 

Using nc Hankel matrices, we can obtain the dual SDP for trace minimization, which is
\begin{equation}
	\begin{split}
		L_{\Theta}(f) =& \inf \ \lr{\M_{d}(L), \mat{G}_{f}}, \\
		& \st \ \M_{d}(L)\succeq0, \\
		& \sts \ \M_d(L)_{1,1} = 1, \\
		& \sts \ \M_d(L)_{u,v} = \M_d(L)_{r,s} \text{ for all } \cyc{u^{\invol}v}{r^{\invol}s}.
	\end{split}
\hspace{10000pt minus 1fil}
\end{equation}

\subsection{No Duality Gap}
For eigenvalue and trace optimization, we can use either the primal or dual SDP, since in both cases the duality gap $\xi-\psi$ is known to be zero.

\begin{theorem}[\protect{\cite[Theorem 4.1]{burgdorf2016optimization}}]
	The SDP pairs \eqref{eq-SOHShierarchy}-\eqref{eq-evSDPdual}, and \eqref{eq-trhierarchy}-\eqref{eq-trdual} satisfy \emph{strong duality}, namely, $\Lambda_{\min}(f) = \lambda_{\min}(f)$, and $\tr_{\Theta}(f)=L_{\Theta}(f)$.
\end{theorem}


In practice most state-of-the-art SDP solvers use primal-dual methods, relying on both formulations. Moreover, guaranteed convergence to an $\epsilon$-optimal solution (for $\epsilon>0$) rely directly on the presence of strong duality \cite{vandenberghe1996semidefinite}. 

\section{Constraints}\label{S-Constraints}

So far we have considered global (eigenvalue or trace) optimization of nc polynomials, where we consider all $\mat{\ul{A}}\in\mbb{S}^{n}$. However, for some applications (such as violation of Bell inequalities, see \S \ref{S-app}) we want to optimize over a specific region of $\mbb{S}^{n}$. 

We now show how to optimize nc polynomials over regions of $\mbb{S}^{n}$ which can be defined via nc polynomial inequalities.  

\subsection{Algebraic extensions}\label{SS-AG}

Let $\mf{g}=\set{ g_{1}, \dotsc, g_{m}}$ be a subset of $\SX$. The \emph{semialgebraic set} associated to $\mf{g}$ is defined as
$$
\mc{D}_{\mf{g}} = \set{ \mat{\ul{A}}\in\mbb{S}^{n} : \forall g\in \mf{g},  g(\mat{\ul{A}})\succeq0 }
$$
i.e., the set of all matrix tuples $\mat{\ul{A}}\in\mbb{S}^{n}$ such that each nc polynomial in $\mf{g}$ has a psd evaluation at $\mat{\ul{A}}$. 

We can naturally extend this notion from matrix tuples of the same order, to bounded self-adjoint operators on some Hilbert space $\mc{H}$, which make $g(\mat{\ul{A}})$ psd for all $g\in \mf{g}$. 
This extension is called the \emph{operator semialgebraic set} associated to  $\mf{g}$, and we denote it by $\mc{D}_{\mf{g}}^{\infty}$.

The \emph{quadratic module} generated by $\mf{g}$ is the set
$$
Q(\mf{g}) := \set{ 
	\sum_{i=0}^{M}  p_{i}^{\invol} g_{i} p_{i} 								\;\middle\vert\; 
		M\in\mbb{N}, g_{i}\in\mf{g}\cup\set{1}, p_{i}\in\RX
}.
$$
Note that when $\mf{g}=\emptyset$, $Q(\emptyset)=\Sigma$. Given $k\in\mbb{N}$, the $k^{\text{th}}$-order truncation of $Q(\mf{g})$ is the following
$$
Q(\mf{g})_{k} := \set{ \sum_{i=0}^{M}  p_{i}^{\invol} g_{i} p_{i} 
				\;\middle\vert\; 
					\begin{aligned}
						M\in\mbb{N}, g_{i}\in\mf{g}\cup\set{1}, 
						p_{i}\in\RX, \\
						\deg(p_{i}^{\invol}g_{i}p_{i})\leq 2k 
					\end{aligned}
	}
$$

We say that $Q(\mf{g})$ is \emph{Archimedean} if for all $q\in\RX$, there is a positive $R\in\mbb{N}$ such that $R-q^{\invol}q\in Q(\mf{g})$. This is equivalent to the existence of an $R\in\mbb{N}$ such that $R-(x_{1}^{2}+\dotsb +x_{n}^{2})\in Q(\mf{g})$ \cite[Proposition 2.2]{klep2007nichtnegativstellensatz}.

The next proposition is obvious, but crucial to the study of nonnegative nc polynomials.
\begin{proposition}[\protect{\cite[Proposition 1.25]{burgdorf2016optimization}}]\label{P-positivityset}
	Let $f\in\SX$ and $\mf{g}\subseteq\SX$. If $f\in Q(\mf{g})$, then $f(\mat{\ul{A}})\succeq0$ for all $\mat{\ul{A}}\in\mc{D}_{\mf{g}}$. Likewise, $f(\mat{\ul{A}})\succeq0$ for all $\mat{\ul{A}}\in\mc{D}_{\mf{g}}^{\infty}$.
\end{proposition}
In general, the converse of Proposition \ref{P-positivityset} is false. However, if we restrict to $\mc{D}_{\mf{g}}^{\infty}$, then we have the noncommutative version of Putinar's Positivstellensatz \cite{helton2004positivstellensatz}.

\begin{theorem}[\protect{\cite[Theorem 1.32]{burgdorf2016optimization}}]\label{T-PutinarPositivstellensatz}
	Let us consider $f\in\SX$ and $\mf{g}\subseteq\SX$ and suppose $Q(\mf{g})$ is Archimedean. If $f\succ0$ for all $\mat{\ul{A}}\in\mc{D}_{\mf{g}}^{\infty}$, then $f\in Q(\mf{g})$.
\end{theorem}

Proposition \ref{P-positivityset} and Theorem \ref{T-PutinarPositivstellensatz} make it possible to develop an approximation hierarchy for constrained nc optimization. 
Before presenting this hierarchy, let us clarify some minor issues that arise with trace optimization.

First, to handle cyclic equivalences in trace optimization, we must use the \emph{cyclic} quadratic module
$$
Q^{\CYC}(\mf{g}) := \set{ 
	f\in\SX 
	\;\middle\vert\;
	\exists g\in Q(\mf{g}) \text{ with } \cyc{f}{g}	
}
$$
and its truncation of order $k$, $Q^{\CYC}(\mf{g})_{k}$, which is defined similarly to $Q(\mf{g})_{k}$.

Second, and perhaps more importantly, we cannot in general consider trace optimization over the operator semialgebraic set $\mc{D}_{\mf{g}}^{\infty}$. This is because, if the space $\mc{H}$ is not finite dimensional, then the algebra $\mc{B}(\mc{H})$ of bounded operators on $\mc{H}$ does not admit a trace. So instead we consider subalgebras of $\mc{B}(\mc{H})$ which admit a trace, namely, finite von Neumann algebras of type $\I$ and $\II$. Then $\mc{D}_{\mf{g}}$ is generalized to the \emph{von Neumann semialgebraic set} $\mc{D}_{\mf{g}}^{\II_{1}}$ (for details on this, we refer the reader to \cite[Definition 1.59]{burgdorf2016optimization} and \cite[Chapter 5]{takesaki2013theory}). 

With this in mind, we have analogous results to Proposition \ref{P-positivityset} and Theorem \ref{T-PutinarPositivstellensatz} for trace positivity. 
\begin{proposition}[\protect{\cite[Proposition 1.62]{burgdorf2016optimization}}]\label{P-trpositivityset}
	Let $f\in\SX$ and $\mf{g}\subseteq\SX$. If $f\in Q^{\CYC}(\mf{g})$, then $\tr(f(\mat{\ul{A}}))\geq0$ for all $\mat{\ul{A}}\in\mc{D}_{\mf{g}}$. Likewise, $\tr(f(\mat{\ul{A}}))\geq0$ for all $\mat{\ul{A}}\in\mc{D}_{\mf{g}}^{\II_{1}}$.
\end{proposition}

\begin{theorem}[\protect{\cite[Proposition 1.63]{burgdorf2016optimization}}]\label{T-tracialPositivstellensatz}
	Let us consider $f\in\SX$ and $\mf{g}\subseteq\SX$ and suppose $Q(\mf{g})$ is Archimedean. 
	If $\tr(f(\underline{A})) > 0$ for all $\underline{A}\in\mc{D}_{\mf{g}}^{\II_{1}}$, then $f \in Q^{\CYC}(\mf{g})$.
\end{theorem}

\subsection{Approximation Hierarchies}\label{SS-approxHierachy}

Given a set of constraints $\mf{g}=\set{g_{1}, \dotsc, g_{m}}$, using Proposition \ref{P-positivityset} we can obtain the following hierarchy for constrained eigenvalue optimization 
\begin{align}\label{eq-consSOHShierarchy} 
	\begin{split}
		\lambda_{\min}^{(d)}(f) = & \sup \ \lambda,\\
		&\st \ f-\lambda\in Q(\mf{g})_{2d}.
	\end{split}
\end{align}

For $g \in \mf{g}$, let us define $d_g := \lceil \deg (g)/2 \rceil$.
To obtain the dual program, we would need the representation of the dual cone $Q(\mf{g})^{\dual}_{2d}$. This dual is nothing more than symmetric linear functionals which are nonnegative on $Q(\mf{g})_{2d}$; in other words, symmetric linear functionals $L$ such that for all $g\in\mf{g}\cup\set{1}$ and $p\in\RX_{d-d_g}$
$$
L( p^{\invol}gp)\geq0.
$$
Similar to the nc Hankel matrix, we define the \emph{nc localizing matrix} $\M_{d-d_g}(gL)$ as 
$$
\M_{d - d_g}(gL)_{u,v} = L(u^{\invol}gv)
$$
(indexed by words $u,v\in\LXR_{d-d_g}$). 
With this in place, we obtain the dual SDP to \eqref{eq-consSOHShierarchy},
\begin{equation}\label{eq-consevSDPdual}
	\begin{split}
		\Lambda_{\min}^{(d)}(f) = & \inf \ \lr{\M_d(L), \mat{G}_{f}}, \\
		& \st \ \M_d(L)_{1,1} = 1, \\
		& \sts  \ \forall g\in\mf{g}\cup\set{1}, \M_{d - d_g}(gL)\succeq0, \\
		& \sts \ \M_d(L)_{u,v} = \M_d(L)_{r,s} \text{ for all } u^{\invol}v = r^{\invol}s.
	\end{split}
\hspace{10000pt minus 1fil}
\end{equation}

With similar reasoning, we obtain the following primal program for trace optimization
\begin{equation*}
	\begin{split}
		\tr_{\Theta}^{(d)} = & \sup \ \tau, \\
		& \st \ f-\tau \in Q^{\CYC}(\mf{g})_{2d},
	\end{split}
\end{equation*}
and the corresponding dual program
\begin{equation*}
	\begin{split}
		L_{\Theta}^{(d)}(f) = & \inf \ \lr{\M_d(L), \mat{G}_{f}}, \\
		& \st \  \M_d(L)_{1,1} = 1,\\
		& \sts \ \forall g\in\mf{g}\cup\set{1}, \M_{d-d_g}(gL)\succeq0,\\
		& \sts \ \M_d(L)_{u,v} = \M_d(L)_{r,s} \text{ for all } \cyc{u^{\invol}v}{r^{\invol}s}.
	\end{split}
	\hspace{10000pt minus 1fil}
\end{equation*}
Furthermore if $Q(\mf{g})$ is Archimedean, then both of these programs converge to the true solution as $d\rightarrow\infty$.

\section{Applications}\label{S-app}

\subsection{Maximal Bell violation}

Consider a quantum mechanical state of two particles, which interacted at some point in the past and are then placed at two separate (far away) locations $A$ and $B$. 
At these locations, we perform some possible measurements $m_{A}$ and $m_{B}$, and get one of two possible outputs $o_{A}$ and $o_{B}$. Because the particles have interacted in the past, the outcomes of the measurements can be correlated, and this correlation can be described by a joint probability $P(m_{A}, m_{B}|o_{A}, o_{B})$. The \emph{generalized Bell inequalities} are inequalities of the form $\sum a_{mo}P(m_{A}, m_{B}|o_{A}, o_{B}) \leq C$, and provide a bound on the possible correlations; which are respected in the setting of classical mechanics, but not in quantum mechanics.

Violation of Bell inequalities serves as an indicator for \emph{entanglement} of a quantum state, which can have many interesting applications such as quantum teleportation \cite{bouwmeester1997experimental}, quantum computation \cite{jozsa2003role} and quantum cryptography \cite{ekert1991quantum}. Moreover, there is interest in understanding the maximal violation of a Bell inequality, for instance \cite{acin2007device} uses this maximum violation to bound the Holevo information in quantum cryptography. 

Since the measurements can be represented as Hermitian operators acting on the quantum state $\mbf{v}$, the expectation of the joint probabilities can be written in the form $E[P(m_{A}, m_{B}|o_{A}, o_{B})] = \mbf{v}^{T}x_{i}y_{j}\mbf{v}$ (where $i,j$ correspond to the possible outputs), and hence we can frame Bell inequality violation as an NCPOP. The most famous Bell inequality is the CHSH inequality \cite{clauser1969proposed}
$$
x_{1}y_{1} + x_{1}y_{2} + x_{2}y_{1} - x_{2}y_{2}\leq 2
$$ 
under the constraints $x_{i}^{2} = y_{j}^{2} = 1$ and $x_{i}y_{j} = y_{j}x_{i}$ for $i,j\in\set{1,2}$. Letting $g(x,y)=x_{1}y_{1} + x_{1}y_{2} + x_{2}y_{1} - x_{2}y_{2}$, 
we can solve the NCPOP
\begin{equation*}
	\begin{split}
		& \inf \ -(g(x,y)+g(x,y)^{\invol})/2 , \\
		& \st \ x^{2}_{i}=1,  \\
		& \sts \ y_{j}^{2} = 1, \\
		& \sts \ x_{i}y_{j}=y_{j}x_{i} \text{ for all } i,j\in\set{1,2},
	\end{split}
\end{equation*}
by first reformulating it to the SOHS program \eqref{eq-consSOHShierarchy} or its dual \eqref{eq-consevSDPdual} as in \S \ref{SS-approxHierachy}, and then solving the associated SDP. The package NCSOStools allows us to do this conveniently with the command \texttt{NCeigMin}, and using a relaxation order of $d=2$, we find that the maximum possible violation of the CHSH inequality is $2\sqrt{2}$.

\subsection{Matrix factorization ranks}

Let $\mat{M}$ be an entrywise nonnegative matrix of size $(p,q)$. Do there exist matrices $\mat{A}_{1},\dotsc,\mat{A}_{p}, \mat{B}_{1}, \dotsc, \mat{B}_{q}\in\mbb{S}_{r}^{\succeq}$ such that we can write $\mat{M}_{ij} = \tr(\mat{A}_{i}\mat{B}_{j})$, and if so, what is the smallest $r$ which allows this? Such a factorization is known as a \emph{psd factorization}, and the smallest admissible $r$ is known as the \emph{psd rank} of $\mat{M}$, denoted  $\rank_{psd}(\mat{M})$ \cite{fiorini2012linear}.

In the context of optimization, the psd rank can be used to examine the size of the \emph{psd lift} of a polytope \cite{thomas2018spectrahedral,gouveia2013lifts}. Letting $\mc{L}$ be an affine subspace of $\mbb{S}_{r}$, the set $\mc{L}\cap\mbb{S}_{r}^{\succeq}$ is known as a \emph{spectrahedron}, which can be written as $\mc{L}\cap\mbb{S}_{r}^{\succeq}=\set{ \mbf{y}\in\mbb{R}^{d} \;\middle\vert\; \mat{A}_{0}+y_{1}\mat{A}_{1}+\dotsb+y_{d}\mat{A}_{d}\succeq0 }$ for some $\mat{A}_{0},\dotsc,\mat{A}_{d}\in\mbb{S}_{r}$. If a polytope $P$ can be written as the projection of a spectrahedron, i.e., $P=\pi(\mc{L}\cap\mbb{S}_{r}^{\succeq})$, then $\mc{L}\cap\mbb{S}_{r}^{\succeq}$ is a psd lift of $P$ of size $r$.

Minimizing a linear functional $L$ over $P$ is the same as minimizing $L\circ \pi$ over $\mc{L}\cap\mbb{S}_{r}^{\succeq}$ (a semidefinite program). If $P$ is defined by $k$ hyperplanes, and has a psd lift $\mc{L}\subseteq\mbb{S}_{r}$ of smaller size ($r\ll k$), then the semidefinite program from the psd lift can be more efficient than the original linear program; this is because the runtime of many linear programming algorithms depends on $k$. Therefore, given a polytope $P$, it is important to know what the lowest possible value for $r$ is, and this can be conveniently described via the psd rank of an associated \emph{slack matrix} $\mat{S}_{P}$ (see \cite[Defintion 3.5]{gouveia2015worst}). 

The psd rank has been broadly studied in many works such as \cite{fawzi2015positive,gouveia2013polytopes,gribling2017matrices,gouveia2015worst}. In \cite{gribling2019lower} the authors present a framework of computing lower bounds for the psd rank (among other factorization ranks), as a hierarchy of NCPOPs. Given a matrix $\mat{M}$ of size $(p,q)$, in our notation this framework can be written as follows, in the nc letters $x_{1},\dotsc,x_{p+q}$. 

Let $\mf{g}$ be the following set in $\SX$,
\begin{align*}
\mf{g}= & \set{ x_{i} - x_{i}^{2} \;\middle\vert\; i=1,\dotsc,p} \\ &\bigcup \set{ \left( \sum_{i=1}^{p} \mat{M}_{ij} \right)x_{p+j} - x_{p+j}^{2}  \;\middle\vert\; j=1,\dotsc,q },
\end{align*}
and 
$$
h(\ul{x}) = 1 -\sum_{i=1}^{p}x_{i}.
$$
Then a sequence of approximations $\rho^{(d)}(\mat{M})$ which bound $\rank_{psd}(\mat{M})$ are computed as the following feasibility problem
\begin{equation}\label{eq-rankpsd}
	\begin{split}
		\rho^{(d)}(\mat{M}) = & \inf \ \M_d(L)_{1,1}, \\
		& \st \ \M_d(L)_{x_{i},x_{p+j}} = \mat{M}_{ij},  \\
		& \sts \  \forall g\in\mf{g}\cup\set{1}, \M_{d - d_g}(gL)\succeq0,\\
		& \sts \ \M_{d - d_h}(hL)=0,\\
		& \sts \ \M_d(L)_{u,v} = \M_d(L)_{r,s} \text{ for all } \cyc{u^{\invol}v}{r^{\invol}s}.
	\end{split}
	\hspace{10000pt minus 1fil}
\end{equation}
\begin{example}
	Consider the following matrix 
	$$
	\mat{M} = 
	\begin{bmatrix}
		1 & 1.75 & 0\\
		0 & 1 & 1.75\\
		1.75 & 0 & 1
	\end{bmatrix}.
	$$
	It is shown in \cite{fawzi2015positive} that $\rank_{psd}(\mat{M}) = 3$. Using the SOHS program \eqref{eq-rankpsd}, we obtain the approximations 
	\begin{equation*}
		\begin{gathered}
			\rho^{(2)}(\mat{M}) \approx 1.90903,\\
			\rho^{(3)}(\mat{M}) \approx 1.90903.\\
		\end{gathered}
	\end{equation*}
	While higher order approximations ($\rho^{(4)}(M), 
	\dotso$) may differ, it is important to keep in mind that although the approximations $\rho^{(d)}(\mat{M})$ may converge, they do not necessarily converge to $\rank_{psd}(\mat{M})$.
\end{example}

As we see from these examples, NCPOPs add value to the broader area of optimization, and indeed quantitative science in general.

\section{Further reading}\label{S-further-reading}

Noncommutative optimization is a highly active area of research which draws inspiration from many different disciplines. NCPOPs can be readily used in many branches of engineering and quantum sciences. We have presented in this article only a fundamental introduction to the area and there are many different avenues one can explore from here. 

As we have seen above, SOHS programming is decidedly dependant on the Positivstellens\"{a}tze from free real algebraic geometry. 
A good place to start learning about free real algebraic geometry and the various nc Positivstellens\"{a}tze is \cite{helton2013free}, which offers insights into important topics such as noncommutative convexity, noncommutative rational functions, and noncommutative spectrahedra (matrix solution sets of linear matrix inequalities). Linear matrix inequalities and noncommutative spectrahedra are of particular importance to engineering and control theory, and Positivstellens\"{a}tze related to these are discussed further in \cite{helton2009convex} and \cite{helton2012convexity}.

The latest Positivstellensatz benefiting NCPOPs focuses on \emph{trace polynomials}; a generalized framework combining what we have presented, polynomials in nc letters and their traces, e.g., $x_{1}^{2}\tr(x_{1}x_{2})+x_{2}x_{1}\tr(x_{3})^{2}$. Trace polynomials are of interest in many areas such as quantum information \cite{navascues2008convergent} and free probability \cite{guionnet2014free}. In \cite{klep2020optimization} the authors present novel SOHS representation results for trace polynomials, and use them to build a converging hierarchy of SDP relaxations, which can be used to optimize over trace polynomials. 
The hierarchy can be practically used in the context of quantum information for finding upper bounds on quantum violations of polynomial Bell inequalities 
\cite{pozsgay2017covariance} and for characterizing entanglement of Werner states \cite{huber2021dimension}.
They also show how to extract minimizers in certain cases. The latest results on trace polynomials in \cite{klep2021positive} provide rational SOHS certificates in the univariate setting, analogous to Artin's solution to Hilbert's 17$^{th}$ problem.

SOHS representations allow convenient modelling of NCPOPs, but equally important is the computational task of solving the SDP relaxations which arise. If the number of variables or relaxation degrees grow too large, programs such as \eqref{eq-evSDPdual} become unmanageable. The \emph{Newton Chip method} \cite{klep2010semidefinite} ameliorates this to some extent, by removing monomials which cannot occur in any SOHS decomposition for a given objective nc polynomial. A cyclic extension of the Newton Chip method for use in trace optimization is given in \cite{burgdorf2013algorithmic}.

The very recent works \cite{klep2021sparse} and \cite{nctssos} analyze the related concept of \emph{correlative and term sparsity} and report significant 
improvements in SDP programs when this is exploited. 
Given an NCPOP in a large number of variables, the underlying idea in this approach is to partition the objective and constraint polynomials in a suitable way, so each partition only depends on a small subset of variables. 
A specific library is available within the TSSOS software \cite{magron2021tssos}.

With similar motivations, the work of \cite{mai2021constant} show that the SDP variables in NCPOPs with a ball constraint, have a constant trace. This allows the use of first order spectral methods to solve the generated SDP, and looks to be a promising approach to reducing the associated computational expense.




\end{document}